\title{\vspace{-0.6cm} The size of a hypergraph and its matching number}
\author{
Hao Huang\thanks{Department of Mathematics, UCLA, Los
Angeles, CA 90095. Email: {\tt huanghao@math.ucla.edu}. Research supported by a
UC Dissertation Year Fellowship.}
\and
Po-Shen Loh\thanks{Department of Mathematical Sciences, Carnegie Mellon University, Pittsburgh, PA 15213.
Email: {\tt ploh@cmu.edu}. Research supported by an NSA Young Investigators Grant.}
\and
Benny Sudakov\thanks{Department of Mathematics, UCLA, Los Angeles, CA 90095. Email:
{\tt bsudakov@math.ucla.edu}. Research supported in part by NSF grant DMS-1101185,
NSF CAREER award DMS-0812005, and by a USA-Israeli BSF grant.}}

\documentclass[11pt]{article}
\usepackage{amsfonts,amssymb,amsmath,latexsym}
\usepackage{verbatim}

\date{}

\newtheorem{thm}{Theorem}[section]

\newtheorem{lemma}[thm]{Lemma}
\newtheorem{cor}[thm]{Corollary}

\newtheorem{conj}[thm]{Conjecture}
\newtheorem{ques}[thm]{Question}

\newenvironment{pf}
      {\medskip\noindent{\bf Proof.}\hspace{1mm}}
      {\hfill$\Box$\medskip}
\newenvironment{proofof}[1]
      {\medskip\noindent{\bf Proof of #1.}\hspace{1mm}}
      {\hfill$\Box$\medskip}
\def\qed{\ifvmode\mbox{ }\else\unskip\fi\hskip 1em plus 10fill$\Box$}

\begin{document}
\maketitle

\begin{abstract}
More than forty years ago, Erd\H{o}s conjectured that for any
$t \leq \frac{n}{k}$, every $k$-uniform hypergraph on $n$ vertices
without $t$ disjoint edges has at most $\max\big\{\binom{kt-1}{k},
\binom{n}{k}-\binom{n-t+1}{k}\big\}$ edges.  Although this appears
to be a basic instance of the hypergraph Tur\'an problem (with a
$t$-edge matching as the excluded hypergraph), progress on this question has remained
elusive.  In this paper, we verify this conjecture for all $t
< \frac{n}{3k^2}$.  This improves upon the best previously known
range $t = O\big(\frac{n}{k^3}\big)$, which dates back to the 1970's.
\end{abstract}

\section{Introduction}
\label{section_introduction}
A $k$-uniform hypergraph is a pair $H=(V,E)$, where $V=V(H)$ is a
finite set of vertices, and $E=E(H)\subseteq \binom Vk$ is a family of
$k$-element subsets of $V$ called edges. A matching in $H$ is a set of
disjoint edges in $E(H)$. We denote by $\nu(H)$ the size of the
largest matching, i.e., the maximum number of disjoint edges in
$H$. The problem of finding the maximum matching in a hypergraph has
many applications in various different areas of mathematics, computer
science, and even computational chemistry.  Yet although the graph
matching problem is fairly well-understood, and solvable in polynomial
time, most of the problems related to hypergraph matching tend to be
very difficult and remain unsolved.  Indeed, the hypergraph matching
problem is known to be NP-hard even for $3$-uniform hypergraphs, without any good approximation
algorithm.

One of the most basic open questions in this area
was raised in 1965 by Erd\H{o}s \cite{erdos-matching}, who asked to
determine the maximum possible number of edges that can appear in any
$k$-uniform hypergraph with matching number $\nu(H) < t \leq
\frac{n}{k}$ (equivalently, without any $t$ pairwise disjoint
edges). He conjectured that this problem has only two extremal
constructions. The first one is a clique consisting of all the $k$-subsets on $kt-1$ vertices,
which obviously has matching number $t-1$. The second example
is a $k$-uniform hypergraph on $n$ vertices containing all the edges intersecting a fixed set of $t-1$ vertices,
which also forces the matching number to be at most
$t-1$.  Neither construction is uniformly better than the other across
the entire parameter space, so the conjectured bound is the maximum
of these two possibilities. Note that in the second case, the complement of this hypergraph
is a clique on $n-t+1$ vertices together with $t-1$ isolated vertices, and thus the original hypergraph has $\binom{n}{k}- \binom{n-t+1}{k}$ edges.

\begin{conj} \label{ec} Every $k$-uniform hypergraph $H$ on $n$
  vertices with matching number $\nu(H)<t \leq \frac{n}{k}$ satisfies
\begin{equation}
  \label{eci}
  e(H)
  \leq
  \max \left\{
    \binom{kt-1}{k} \,,
    \binom{n}{k}-\binom{n-t+1}{k}
  \right\} \,.
\end{equation}
\end{conj}

In addition to being important in its own right, this Erd\H{o}s conjecture has several interesting applications, which we discuss in the concluding remarks.
Yet although it is more than forty years old, only partial results have been discovered so far.  In the case
$t=2$, the condition simplifies to the requirement that every pair of
edges intersects, so Conjecture \ref{ec} is thus equivalent to a
classical theorem of Erd\H{o}s, Ko, and Rado \cite{erdos-ko-rado}:
that any intersecting family of $k$-subsets on $n \geq 2k$ elements
has size at most $\binom{n-1}{k-1}$.  The graph case ($k=2$) was
separately verified in \cite{erdos-gallai} by Erd\H{o}s and Gallai.
For general fixed $t$ and $k$, Erd\H{o}s \cite{erdos-matching} proved
his conjecture for sufficiently large $n$. Frankl \cite{frankl-survey} showed that Conjecture \ref{ec}
was asymptotically true for all $n$ by proving the weaker bound $e(H) \leq (t-1) \binom{n-1}{k-1}.$

A short calculation shows that when $t \leq \frac{n}{k+1}$, we
always have $\binom{n}{k} - \binom{n-t+1}{k} > \binom{kt-1}{k}$, so
the potential extremal example in this case has all edges intersecting a fixed set of $t-1$ vertices.
One natural question is then to determine the range of $t$ (with respect
to $n$ and $k \geq 3$) for which the maximum is indeed equal to
$\binom{n}{k}-\binom{n-t+1}{k}$, i.e., where the second case is
optimal. Recently, Frankl, R\"{o}dl, and Ruci\'{n}ski
\cite{frankl-rodl-rucinski} studied 3-uniform hypergraphs ($k=3$), and
proved that for $t \leq n/4$, the maximum was indeed $\binom{n}{3} -
\binom{n-t+1}{3}$, establishing the conjecture in that range.
For general $k \geq 4$, Bollob\'{a}s, Daykin, and Erd\H{o}s
\cite{bollobas-daykin-erdos} explicitly computed the bounds achieved
by the proof in \cite{erdos-matching}, showing that the conjecture
holds for $t<\frac{n}{2k^3}$.  Frankl and Fur\"{e}di
\cite{frankl-survey} established the result in a different range $t <
\big(\frac{n}{100k}\big)^{1/2}$, which improves the original bound
when $k$ is large relative to $n$.  In this paper, we
extend the range in which the Erd\H os conjecture holds to all $t<\frac{n}{3k^2}$.

\begin{thm}\label{main}
  For any integers $n, k, t$ satisfying $t<\frac{n}{3k^2}$, every
  $k$-uniform hypergraph on $n$ vertices without $t$ disjoint edges
  contains at most $\binom{n}{k}-\binom{n-t+1}{k}$ edges.
\end{thm}

To describe the idea of our proof, we first outline Erd\H{o}s's
original approach for the case $t < \frac{n}{2k^3}$. Let $v$ be a
vertex of maximum degree.  By induction on $t$ we find $t-1$ disjoint
edges $F_1, \ldots, F_{t-1}$, none of which contain $v$. If $\deg(v)$
exceeds $k(t-1)\binom{n-2}{k-2}$, which is the maximum possible number
of edges containing $v$ which also meet a vertex in
$\bigcup_{i=1}^{t-1}F_i$, then we can find $t$ disjoint
edges. Otherwise, the number of edges meeting any of $F_i$ is at most
$|\bigcup_{i=1}^{t-1}F_i| \cdot k(t-1)\binom{n-2}{k-2}=k(t-1)\cdot
k(t-1)\binom{n-2}{k-2}$, which turns out to be less than the total
number of edges when $n \geq 2k^3 t$.  Any other edge will serve as
the $t$-th edge in the matching.

To improve Erd\H{o}s's bound, we show that in the first part of the
argument, we are already done if the $t$-th largest degree exceeds $2t
\binom{n-2}{k-2}$.  This puts a tighter constraint on the sum of the
degrees of the $k(t-1)$ vertices in $\bigcup_{i=1}^{t-1}F_i$, allowing
the second stage to proceed under the relaxed assumption $n \geq
3k^2t$.  The fact that $t$ vertices of degree at least $2t
\binom{n-2}{k-2}$ are enough to find $t$ disjoint edges leads
naturally to the following multicolored version of the Erd\H os
conjecture, which was also considered independently by Aharoni and
Howard in \cite{aharoni}.

\begin{conj}\label{ahc}
Let $\mathcal{F}_1, \ldots, \mathcal{F}_t$ be families of subsets in $\binom{[n]}{k}$. If $|\mathcal{F}_i|>
\max \left \{\binom{n}{k}-\binom{n-t+1}{k}, \binom{kt-1}{k}\right\}$ for all $1 \leq i\leq t$, then there is a ``rainbow'' matching of size $t$: one that contains
exactly one edge from each family.
\end{conj}
The $k=2$ case of this conjecture was established by Meshulam (see \cite{aharoni}).
To obtain Theorem \ref{main}, we prove an asymptotic version of Conjecture \ref{ahc}, by showing that a rainbow matching exists whenever $|\mathcal{F}_i| > (t-1) \binom{n-1}{k-1}$ for every $1
\leq i \leq t$.

The rest of this paper is organized as follows. In the next section,
we describe the so-called  \emph{shifting} method, which is a well known technique in extremal set theory, and use it to prove some
preliminary results.  In Section \ref{section_maintheorem} we
first prove the multicolored
Erd\H{o}s conjecture asymptotically, and then use it to prove Theorem \ref{main}.
There, we also use the same argument to show that Conjecture \ref{ahc} holds for all $t < \frac{n}{3k^2}$.
The last section contains some concluding remarks and open problems.

\section{Shifting}
\label{section_shift}

In extremal set theory, one of the most important and widely-used
tools is the technique of shifting,
which allows us to limit our attention to sets with certain structure.
In this section we will only state and prove the relevant results for
Section \ref{section_maintheorem}.  For more background on the
applications of shifting in extremal set theory, we refer the reader
to the survey \cite{frankl-survey} by Frankl.

Given a family $\mathcal{F}$ of equal-size subsets of $[n]$, for
integers $1 \leq j < i \leq n$, we define the $(i, j)$-shift map
$S_{ij}$ as follows: for any set $F \in \mathcal{F}$,
\begin{displaymath}
S_{ij}(F)=
\begin{cases}
F \setminus \{i\} \cup \{j\} \,,& \textrm{iff~} i \in F,~j \not\in F~\textrm{and~} F \setminus \{i\} \cup \{j\} \not \in \mathcal{F}\,;\\
F\,,&\textrm{otherwise.}
\end{cases}
\end{displaymath}
Also, we denote the family after shifting as
\begin{displaymath}
  S_{ij}(\mathcal{F}) = \left\{
    S_{ij}(F): F \in \mathcal{F}
  \right\} \,.
\end{displaymath}

\begin{lemma}
  \label{lemma_shift}
  The shift map $S_{ij}$ satisfies the following properties.
  \begin{description}
  \item[(i)] $|S_{ij}(\mathcal{F})|=|\mathcal{F}|$.
  \item[(ii)] If $\mathcal{F}$ is $k$-uniform, then so is
    $S_{ij}(\mathcal{F})$.
  \item[(iii)] If the families $\mathcal{F}_1$, \ldots,
    $\mathcal{F}_t$ have the property that no subsets $F_1 \in
    \mathcal{F}_1$, \ldots, $F_t \in \mathcal{F}_t$ are pairwise
    disjoint, then the shifted families $S_{ij}(\mathcal{F}_1)$,
    \ldots, $S_{ij}(\mathcal{F}_t)$ still preserve this property.
  \end{description}
\end{lemma}

\begin{pf}
  Claims (i) and (ii) are obvious. For (iii), assume that the
  statement is false, i.e., we have $F_i \in \mathcal{F}_i$ such that
  $S_{ij}(F_1)$, \ldots, $S_{ij}(F_t)$ are pairwise disjoint, while
  $F_1$, \ldots, $F_t$ are not.  Without loss of generality, $F_1 \cap
  F_2 \neq \emptyset$.  Next, observe that whenever $S_{ij}(F_k) \neq
  F_k$, we also have $j \in S_{ij}(F_k)$, so the pairwise disjointness
  of the $S_{ij}(F_k)$ implies that the only possible case
  (re-indexing if necessary) is for $S_{ij}(F_1) = F_1 \setminus
  \{i\} \cup \{j\}$, and $S_{ij}(F_k) = F_k$ for every $k \geq 2$.
  Note also that since $F_1$ and $F_2$ intersect while $S_{ij}(F_1)$
  and $S_{ij}(F_2)$ do not, we must have $i \in F_2$ and $j \not \in
  F_2$.

  Therefore the only reason that $S_{ij}(F_2)=F_2$ is because $F'_2 =
  F_2 \setminus \{i\} \cup \{j\}$ is already in $\mathcal{F}_2$.  The
  pair of disjoint sets $S_{ij}(F_1)$ and $S_{ij}(F_2) = F_2$ have the
  same union as the pair of disjoint sets $F_1$ and $F_2'$.  Using the
  pairwise disjointness of the $S_{ij}(F_k)$, we conclude that the
  sets $F_1$, $F_2'$, $F_3$, \ldots, $F_t$ are pairwise disjoint as
  well, contradicting our initial assumption.
\end{pf}

In practice, we often combine the shifting technique with induction on
the number of elements in the underlying set.  Indeed, let us apply
the shifts $\{S_{ni}\}_{1 \leq i \leq n-1}$ successively, and with
slight abuse of notation, let us again call the resulting families
$\mathcal{F}_1$, \ldots, $\mathcal{F}_t$.  Create from each
$\mathcal{F}_i$ two sub-families based on containment of the final
element $n$:
\begin{align*}
  \mathcal{F}_i(n) &= \left\{
    F \setminus \{n\} : F \in\mathcal{F}_i, n \in F
  \right\} \,, \\
  \mathcal{F}_i(\bar{n}) &= \left\{
    F \hspace{27.7pt}: F \in \mathcal{F}_i, n \not \in F
  \right\} \,.
\end{align*}
It turns out that the rainbow matching number does not increase by this decomposition.

\begin{lemma}
\label{induction_shifting}
Let $\mathcal{F}_1$, \ldots, $\mathcal{F}_t$ be the shifted families,
where each $\mathcal{F}_i$ is $k_i$-uniform and $\sum_{i=1}^t k_i \leq
n$.  Suppose that no subsets $F_1 \in \mathcal{F}_1$, \ldots, $F_t \in
\mathcal{F}_t$ are pairwise disjoint.  Then, for any $0 \leq s \leq
t$, the families $\mathcal{F}_1(n), \ldots, \mathcal{F}_s(n),
\mathcal{F}_{s+1}(\bar{n}), \ldots, \mathcal{F}_t(\bar{n})$ still have
the same property.
\end{lemma}
\begin{pf}
  Assume for the sake of contradiction that there exist pairwise
  disjoint sets $F_1 \in \mathcal{F}_1(n)$, \ldots, $F_s \in
  \mathcal{F}_s(n)$, $F_{s+1} \in \mathcal{F}_{s+1}(\bar{n})$, \ldots,
  $F_t \in \mathcal{F}_t(\bar{n})$. By definition of
  $\mathcal{F}_i(n)$ and $\mathcal{F}_i(\bar{n})$, we know that $F_i
  \cup \{n\} \in \mathcal{F}_i$ for $1 \leq i \leq s$, and $F_i \in
  \mathcal{F}_i$ for $s+1 \leq i \leq t$.  The size of
  $\bigcup_{i=1}^t F_i$ is equal to
  \begin{displaymath}
    \sum_{i=1}^t |F_i|
    =
    \sum_{i=1}^s (k_i-1) + \sum_{i=s+1}^t k_i
    =
    \sum_{i=1}^t k_i - s
    \leq
    n - s \,,
  \end{displaymath}
  so there exist distinct elements $x_1, \ldots, x_s \not \in \bigcup_{i=1}^t
  F_i$. Since $F_i \cup \{n\}$ is invariant under the shift
  $S_{nx_i}$, the set $F_i \cup \{x_i\} = (F_i \cup \{n\}) \setminus
  \{n\} \cup \{x_i\}$ must also be in the family
  $\mathcal{F}_i$. Taking $F'_i=F_i \cup \{x_i\}$ for $1 \leq i \leq
  s$, together with $F_i$ for $s+1 \leq i \leq t$, it is clear that we
  have found pairwise disjoint sets from $\mathcal{F}_i$,
  contradiction.
\end{pf}

\section{Main result}
\label{section_maintheorem}

In this section, we discuss the Erd\H{o}s conjecture and its
multicolored generalizations, and prove the original conjecture for
the range $t<\frac{n}{3k^2}$.  The colored interpretation arises from
considering the collection of families $\mathcal{F}_i$ as a single uniform
hypergraph (possibly with repeated edges) on the vertex set $[n]$,
where each set in $\mathcal{F}_i$ introduces a hyperedge colored in
the $i$-th color.  The following lemma is a multicolored
generalization of Theorem 10.3 in \cite{frankl-survey}, and provides a
sufficient condition for a multicolored hypergraph to contain a
rainbow matching of size $t$.

\begin{lemma}\label{mainlemma}
  Let $\mathcal{F}_1$, \ldots, $\mathcal{F}_t$ be families of subsets
  of $[n]$ such that for each $i$, $\mathcal{F}_i$ only contains sets
  of size $k_i$, $|\mathcal{F}_i| > (t-1) \binom{n-1}{k_i-1}$, and $n
  \geq \sum_{i=1}^t k_i$.  Then there exist $t$ pairwise disjoint sets
  $F_1 \in \mathcal{F}_1$, \ldots, $F_t \in \mathcal{F}_t$.
\end{lemma}
\begin{pf}
  We proceed by induction on $t$ and $n$.  The case $t=1$ is
  trivial.  For general $t$, we can also handle all minimal cases of
  the form $n= \sum_{i=1}^t k_i$.  Indeed, consider a uniformly random
  permutation $\pi$ of $[n]$, and define a series of indicator random
  variables $\{X_i\}$ as follows: $X_1=1$ iff $\{\pi(1), \ldots,
  \pi(k_1)\}$ is a set in $\mathcal{F}_1$ and $X_1=0$ otherwise, and in general,
 $X_j=1$  iff $\{\pi(k_1+\cdots+k_{j-1}+1), \ldots, \pi(k_1+\cdots+k_j)\}$ is a set in
  $\mathcal{F}_j$. We assume that there are no $t$ disjoint sets
  from different families, so we deterministically have:
  \begin{equation}\label{matchingbound}
    X_1 + \cdots + X_t \leq t-1 \,.
  \end{equation}
  On the other hand, it is easy to see that
the expectation of $X_i$ is the probability that
  a random $k_i$-set is in $\mathcal{F}_i$, so
  \begin{displaymath}
    \mathbb{E} X_i= \dfrac{|\mathcal{F}_i|}{\binom{n}{k_i}} \,.
  \end{displaymath}
  Yet we know that for every $i$, we have $|\mathcal{F}_i|>(t-1)
  \binom{n-1}{k_i-1}$, so
  \begin{displaymath}
    \mathbb{E} X_i
    >
    \dfrac{(t-1)\binom{n-1}{k_i-1}}{\binom{n}{k_i}}
    =
    (t-1) \frac{k_i}{n} \,.
  \end{displaymath}
  Summing these inequalities over $1 \leq i \leq t$, we obtain that $\sum_{i=1}^t  \mathbb{E} X_i>t-1$, a
  contradiction to \eqref{matchingbound}.

  Now we consider a generic instance with $n > \sum_{i=1}^t k_i$, and
  inductively assume that all instances with smaller $n$ are known.
  By Lemma \ref{lemma_shift}, after applying all shifts $\{S_{ni}\}_{1
    \leq i \leq n-1}$, we obtain families in which any rainbow $t$-matching can be pulled back to
  a rainbow $t$-matching in $\{\mathcal{F}_i\}$. For convenience we still call the shifted families $\{\mathcal{F}_i\}$.
  Our next step
  is to partition each $\mathcal{F}_i$ into $\mathcal{F}_i(n) \cup
  \mathcal{F}_i(\bar{n})$, but in order to avoid empty sets, we first
  dispose of the case when there is some $k_i = 1$ with $\{n\} \in
  \mathcal{F}_i$.  After re-indexing, we may assume that this is
  $\mathcal{F}_1$.  Since $|\mathcal{F}_i|>(t-1) \binom{n-1}{k_i-1}$
  and there are at most $\binom{n-1}{k_i-1}$ sets containing $n$,
  every other $\mathcal{F}_i$ has more than $(t-2) \binom{n-1}{k_i-1}$
  sets which in fact lie in $[n-1]$. By induction on the $t-1$ sizes
  $k_2, \ldots, k_t$, we find $t-1$ such disjoint sets from
  $\mathcal{F}_2$, \ldots, $\mathcal{F}_t$ which, together with $\{n\}
  \in \mathcal{F}_1$, establish the claim.

  Returning to the general case, since
  $|\mathcal{F}_i|=|\mathcal{F}_i(n)|+|\mathcal{F}_i(\bar{n})|$ and
  our size condition is
  \begin{displaymath}
    |\mathcal{F}_i| > (t-1) \binom{n-1}{k_i-1} = (t-1) \binom{n-2}{k_i-2} + (t-1) \binom{n-2}{k_i-1} \,,
  \end{displaymath}
  we conclude that for each $i$, either $|\mathcal{F}_i(n)| > (t-1)
  \binom{n-2}{k_i-2}$ or $|\mathcal{F}_i(\bar{n})| > (t-1)
  \binom{n-2}{k_i-1}$. Without loss of generality, we may assume that
  $|\mathcal{F}_i(n)| > (t-1) \binom{n-2}{k_i-2}$ for $1 \leq i \leq
  s$, and $|\mathcal{F}_i(\bar{n})| > (t-1) \binom{n-2}{k_i-1}$ for
  $s+1 \leq i \leq t$.  Note that $\mathcal{F}_i$ is $(k_i-1)$-uniform
  for $1 \leq i \leq s$ and $k_i$-uniform for $s+1 \leq i \leq t$, and
  the base set now has $n-1$ elements. Induction on $n$ and Lemma
  \ref{induction_shifting} then produce $t$ disjoint sets from
  different families.
\end{pf}

As mentioned in the introduction, the conjectured extremal hypergraph
when $t \leq \frac{n}{k+1}$ is the hypergraph consisting of all edges intersecting a
fixed set of size $t-1$.  If we inspect the vertex degree sequence of this
hypergraph, we observe that although there are $t-1$ vertices with
high degree $\binom{n-1}{k-1}$, the remaining vertices only have
degree $\binom{n-1}{k-1}-\binom{n-t}{k-1}$.  For small $t$, this is
asymptotically about $(t-1) \binom{n-2}{k-2}$, which is much smaller
than $\binom{n-1}{k-1} = \frac{n-1}{k-1} \binom{n-2}{k-2}$. The
following corollary of Lemma \ref{mainlemma} shows that this sort of
phenomenon generally occurs when hypergraphs satisfy the conditions in
the Erd\H{o}s conjecture.

\begin{cor}\label{corollary}
  If a $k$-uniform hypergraph $H$ on $n$ vertices has $t$ distinct
  vertices $v_1$, \ldots, $v_t$ with degrees $d(v_i) > 2(t-1)
  \binom{n-2}{k-2}$, and $kt \leq n$, then $H$ contains $t$ disjoint
  edges.
\end{cor}
\begin{pf}
  Let $H_i$ be a $(k-1)$-uniform hypergraph containing all the subsets of $V(H) \setminus \{v_1, \ldots, v_t\}$ of size $k-1$ which together with
  $v_i$ form an edge of $H$.
For any fixed $1
  \leq i \leq t$ and $j \neq i$, there are at most $\binom{n-2}{k-2}$
  edges of $H$ containing both vertices $v_i$ and $v_j$. Therefore for every
   hypergraph $H_i$,
  \begin{displaymath}
    e(H_i)
    \geq
    d(v_i)- (t-1) \binom{n-2}{k-2}
    >
    (t-1) \binom{n-2}{k-2}
    \geq
    (t-1)\binom{n-t-1}{k-2} \,.
  \end{displaymath}
  Since every hypergraph $H_i$ is $(k-1)$-uniform and has $n-t$
  vertices, we can use Lemma \ref{mainlemma} with $\mathcal{F}_i=E(H_i)$,
  $k_i=k-1$ and  $n$ replaced by $n-t$, to find $t$ disjoint edges
  $e_1 \in E(H_1)$, \ldots, $e_t \in E(H_t)$. Taking the edges $e_i
  \cup \{v_i\} \in E(H)$, we obtain $t$ disjoint edges in the original
  hypergraph $H$.
\end{pf}

Now we are ready to prove our main result, Theorem \ref{main}, which
states that for $t < \frac{n}{3k^2}$, every $k$-uniform hypergraph on
$n$ vertices without $t$ disjoint edges contains at most
$\binom{n}{k}-\binom{n-t+1}{k}$ edges.

\begin{proofof}{Theorem \ref{main}}
  We proceed by induction on $t$.  The base case $t=1$ is trivial, so
  we consider the general case, assuming that the $t-1$ case is known.
  Suppose $e(H) > \binom{n}{k}-\binom{n-t+1}{k}$, and let us seek $t$
  disjoint edges in $H$.  We first consider the situation when there
  is a vertex $v$ of degree $d(v) > k(t-1) \binom{n-2}{k-2}$.  Let
  $H_v$ be the sub-hypergraph induced by the vertex set $V(H)
  \setminus \{v\}$.  Since there are at most $\binom{n-1}{k-1}$ edges
  containing $v$,
  \begin{align*}
    e(H_v)
    \geq
    e(H) - \binom{n-1}{k-1}
    &> \binom{n}{k} - \binom{n-t+1}{k} - \binom{n-1}{k-1} \\
    &= \binom{n-1}{k}-\binom{(n-1)-(t-1)+1}{k} \,.
  \end{align*}
  By induction, there are $t-1$ disjoint edges $e_1$, \ldots,
  $e_{t-1}$ in $H_v$, spanning $(t-1)k$ distinct vertices $u_1$,
  \ldots, $u_{(t-1)k}$. Note that the number of edges containing $v$ and any vertex $u_j$ is
  at most $\binom{n-2}{k-2}$. Therefore since we assumed that
  $d(v) > k(t-1) \binom{n-2}{k-2}$, there must be another edge $e_t$
  which contains $v$ but avoids $u_1$, \ldots, $u_{(t-1)k}$. We then
  have $t$ disjoint edges $e_1, \ldots, e_{t}$ in $H$.

  Now suppose that the maximum vertex degree in $H$ is at most $k(t-1)
  \binom{n-2}{k-2}$.  After re-indexing the vertices, we may assume
  that $k(t-1) \binom{n-2}{k-2} \geq d(v_1) \geq \cdots \geq
  d(v_n)$. If the $t$-th largest degree satisfies $d(v_t) > 2(t-1)
  \binom{n-2}{k-2}$, then Corollary \ref{corollary} immediately
  produces $t$ disjoint edges in $H$, so we may also assume for the
  remainder that $d(v_t) \leq 2(t-1) \binom{n-2}{k-2}$.

  By induction (with room to spare), we also know that there are $t-1$
  disjoint edges in $H$, spanning $(t-1)k$ vertices. Among these
  vertices, the $t-1$ largest degrees are at most $k(t-1)
  \binom{n-2}{k-2}$ by our maximum degree assumption, while the
  remaining $(t-1)(k-1)$ vertices cannot have degrees exceeding
  $d(v_t) \leq 2(t-1) \binom{n-2}{k-2}$.  Therefore the sum of degrees
  of these $(t-1)k$ vertices is at most
  \begin{displaymath}
    (t-1) \cdot k(t-1) \binom{n-2}{k-2} + (t-1)(k-1) \cdot 2(t-1) \binom{n-2}{k-2}
    =
    (t-1)^2(3k-2) \binom{n-2}{k-2} \,.
  \end{displaymath}
  However, we know that the total number of edges exceeds
  \begin{align*}
    e(H) &> \binom{n}{k}-\binom{n-t+1}{k} \\
    &= \left[
      1
      -
      \left( 1-\frac{t-1}{n} \right) \cdots \left( 1-\frac{t-1}{n-k+1} \right)
    \right] \binom{n}{k} \\
    &\geq \left[
      1 - \left( 1-\frac{t-1}{n} \right)^k
    \right] \binom{n}{k} \\
    &\geq \left[
      k \cdot \frac{t-1}{n} - \binom{k}{2} \left( \frac{t-1}{n} \right)^2
    \right] \dfrac{n(n-1)}{k(k-1)} \binom{n-2}{k-2}\\
    &\geq \left(
      \frac{(n-1)(t-1)}{k-1} - \frac{(t-1)^2}{2}
    \right) \binom{n-2}{k-2} \,,
  \end{align*}
where we used that $(1-x)^k \leq 1-kx+{k \choose 2}x^2$ when $0 \leq kx \leq 1$.
  Since $n > 3k^2 t$, we also have $n-1 > 3k(k-1) (t-1)$.  Therefore,
  \begin{displaymath}
    e(H)
    >
    (t-1)^2 \left( 3k - \frac{1}{2} \right) \binom{n-2}{k-2} \,,
  \end{displaymath}
  and so there is another edge in $H$ disjoint from the previous $t-1$
  edges, again producing $t$ disjoint edges in $H$.
\end{proofof}

Based on the same idea and technique, we can also obtain a
multicolored version of the Erd\H{o}s conjecture, which is an
analogue of a theorem of Kleitman \cite{kleitman} for matching number
greater than one.  Note that Theorem \ref{main} is the $\mathcal{F}_1
= \cdots = \mathcal{F}_t$ case of the following result.

\begin{thm}\label{multiversion}
  Let $\mathcal{F}_1$, \ldots, $\mathcal{F}_t$ be $k$-uniform families
  of subsets of $[n]$, where $t < \frac{n}{3k^2}$, and every
  $|\mathcal{F}_i| > \binom{n}{k} - \binom{n-t+1}{k}$.  Then there
  exist pairwise disjoint sets $F_1 \in \mathcal{F}_1$, \ldots, $F_t
  \in \mathcal{F}_t$.
\end{thm}
\begin{pf}
  For any vertex $v \in \mathcal{F}_i$, let $H_v^j$ be the
  sub-hypergraph of $\mathcal{F}_j$ induced by the vertex set $[n]
  \setminus \{v\}$.  Then as in the previous proof,
  \begin{displaymath}
    e(H_v^j) \geq |\mathcal{F}_i|- \binom{n-1}{k-1}
    >
    \binom{n-1}{k}-\binom{(n-1)-(t-1)+1}{k} \,.
  \end{displaymath}
  By induction on $t$, for every $i$ there exist $t-1$ disjoint edges
  $\{e_j\}_{j \neq i}$ such that $e_j \in H_v^j$. So as before, if
  some $\mathcal{F}_i$ has a vertex with degree $d(v) > k(t-1)
  \binom{n-2}{k-2}$, then there is an edge in $\mathcal{F}_i$ which contains $v$ and
  is disjoint from $\{e_j\}_{j \neq i}$. Hence we may assume the
  maximum degree in each hypergraph $\mathcal{F}_i$ is at most $k(t-1)
  \binom{n-2}{k-2}$.

  On the other hand, by induction on $t$ we also know that for every $i$
  there exist $t-1$ disjoint edges from the families
  $\{\mathcal{F}_j\}_{j \neq i}$, spanning $(t-1)k$ vertices. If some $\mathcal{F}_i$ has $t$-th
  largest degree at most $2(t-1)\binom{n-2}{k-2}$, then the sum of
  degrees of these $(t-1)k$ vertices in $\mathcal{F}_i$ is again at
  most
  \begin{displaymath}
    (t-1)^2(3k-2) \binom{n-2}{k-2}
    \leq
    \binom{n}{k}-\binom{n-t+1}{k}
    <
    e(\mathcal{F}_i) \,,
  \end{displaymath}
  which guarantees the existence of an edge in $\mathcal{F}_i$
  disjoint from the previous $t-1$ edges from $\{\mathcal{F}_j\}_{j
    \neq i}$. So, we may assume that each $\mathcal{F}_i$ contains at
  least $t$ vertices with degree above $2(t-1)\binom{n-2}{k-2}$.

  Now select distinct vertices $v_i$, such that for each $1 \leq i \leq t$, the degree of
  $v_i$ in $\mathcal{F}_i$ exceeds $2(t-1) \binom{n-2}{k-2}$.
  Consider all the subsets of $[n] \setminus \{v_1, \ldots, v_t\}$ which together with
$v_i$ form an edge of $\mathcal{F}_i$. Denote this $(k-1)$-uniform hypergraph by $T^i$.
The same calculation as in Corollary \ref{corollary} gives
  \begin{displaymath}
    e(T^i) > (t-1) \binom{n-t-1}{k-2} \,.
  \end{displaymath}
  Applying Lemma \ref{mainlemma} to $\{T^i\}$, we again find $t$
  disjoint edges from different families, as desired.
\end{pf}

\section{Concluding Remarks}
\begin{list}{\labelitemi}{\leftmargin=1em}
\item In this paper, we proved that for $t < \frac{n}{3k^2}$, every
  $k$-uniform hypergraph on $n$ vertices with matching number less
  than $t$ has at most $\binom{n}{k} - \binom{n-t+1}{k}$ edges. This
  verifies the conjecture of Erd\H{o}s in this range of $t$, and
  improves upon the previously best known range by a factor of $k$. As
  we discussed in the introduction, if the Erd\H{o}s conjecture is
  true in general, then for $t < \frac{n}{k+1}$, the maximum number of
  edges cannot exceed $\binom{n}{k} - \binom{n-t+1}{k}$.  It would be
  very interesting to tighten the range to $t< O\big(
  \frac{n}{k}\big)$.

\item \emph{A fractional matching} in a $k$-uniform hypergraph
  $H=(V,E)$ is a function $w:E\to[0,1]$ such that for each $v\in V$ we
  have $\sum_{e\ni v}w(e)\le1$. The \emph{size}\/ of $w$ is the sum
  $\sum_{e\in E}w(e)$, and the size of the largest fractional matching
  in $H$ is denoted by $\nu^*(H)$. The fractional version of the
  Erd\H{o}s conjecture states that among $k$-uniform hypergraphs $H$ on
  $n$ vertices with fractional matching number $\nu^*(H) < xn$, the
  maximum number of edges is asymptotically $(1+o(1))\max
  \big\{(kx)^k, 1-(1-x)^k \big\} \binom{n}{k}$.

  It appears that these conjectures are closely related to several
  other interesting problems. For example, it was shown in
  \cite{six-person} that the integral version can be used to determine
  the minimum degree condition which ensures the existence of perfect
  matchings in uniform hypergraphs. Furthermore, it turns out that the
  fractional version is closely related to an old probability
  conjecture of Samuels \cite{samuels_chebyshev} and in computer
  science, it has applications to finding optimal data allocations in
  distributed storage systems (see \cite{six-person} for more
  details).  In \cite{alon-huang-sudakov}, the fractional Erd\H{o}s
  conjecture was used to attack an old problem of
  Manickam-Mikl\'{o}s-Singhi, which states that for $n \geq 4k$, every
  set of $n$ real numbers with nonnegative sum has at least
  $\binom{n-1}{k-1}$ $k$-element subsets whose sums are also
  nonnegative.

\item Pyber \cite{pyber} proved the following product-type
  generalization of the Erd\H os-Ko-Rado theorem. Let
  $\mathcal{F}_1$ and $\mathcal{F}_2$ be families of $k_1$- and
  $k_2$-element subsets of $[n]$. If every pair of sets $F_1 \in
  \mathcal{F}_1$ and $F_2 \in \mathcal{F}_2$ intersects, then
  $|\mathcal{F}_1||\mathcal{F}_2| \leq \binom{n-1}{k_1-1}
  \binom{n-1}{k_2-1}$ for sufficiently large $n$. The special case
  when $k_1=k_2$ and $\mathcal{F}_1=\mathcal{F}_2$ corresponds to the
  Erd\H{o}s-Ko-Rado theorem.  Our Theorem \ref{multiversion} is a
  minimum-type result of similar flavor.  Hence, it would be interesting to study
 the following multicolor analogue of Pyber's result.
\begin{ques}\label{product_multi}
  What is the maximum of $\prod_{i=1}^t |\mathcal{F}_i|$
  among families $\mathcal{F}_1, \ldots, \mathcal{F}_t$ of subsets
  of $[n]$, where each $\mathcal{F}_i$ is $k_i$-uniform, and
  there are no $t$ pairwise disjoint subsets $F_1 \in \mathcal{F}_1$,
  \ldots, $F_t \in \mathcal{F}_t$?
\end{ques}

\end{list}

\noindent \textbf{Acknowledgments}~~The authors would like to thank the anonymous referee for 
carefully reading the paper and the many helpful comments.

\end{document}